\documentclass[a4paper,9pt,titlepage,leqno,article,oneside,openany]{memoir}
\chapterstyle{article}
\setbeforesubsecskip{1ex}
\setaftersubsecskip{-0.25em}
\setsubsecheadstyle{\normalsize\bfseries}

\usepackage{amsmath}
\usepackage{amsfonts}

\usepackage{txfonts}
\usepackage[numbers]{natbib}
\bibpunct{[}{]}{;}{n}{}{,}
\usepackage{amssymb}
\usepackage{xcolor}
\usepackage{mathptmx}
\usepackage{graphicx}          

\newtheorem{defn}{Definition} 
\newtheorem{lem}[defn]{Lemma} 
\newtheorem{cor}[defn]{Corollary} 
\newtheorem{prop}[defn]{Proposition} 
\newenvironment{mypf}{\emph{Proof.}}{~\hfill\rule{0.75em}{0.75em}\\[-0.0em]}

\newcommand{\eotsymbol}{~\hfill{\raise.5ex\hbox{\fbox{}}}\\[-0.0em]}

\newcommand{\T}{\mathsf{T}}
\newcommand{\tr}[1]{\mathop{\mathrm{tr}}\left({#1}\right)}
\newcommand{\grad}{\mathop{\mathrm{grad}}}
\newcommand{\md}{d}

\title{On the differential equation
$\dot{\Theta}=(\Theta^{\T}-\Theta)\Theta$ with $\Theta\in{SO(n)}$}
\author{Gerd S. Schmidt, Christian Ebenbauer, Frank Allg{\"o}wer\\
Institute for Systems Theory and Automatic Control}

\usepackage{ifthen}
\newcommand{\ForDiss}[1]{\ifthenelse{\boolean{ForDiss}}{{#1}}{}}
\newboolean{ForDiss}
\setboolean{ForDiss}{false}
\newcommand{\ForPaper}[1]{\ifthenelse{\boolean{ForPaper}}{{#1}}{}}
\newboolean{ForPaper}
\setboolean{ForPaper}{true}

\begin{document}
\maketitle
\begin{abstract}
  In this note we consider the global convergence properties of the differential
  equation $\dot{\Theta}=(\Theta^{\T}-\Theta)\Theta$ with $\Theta\in{SO(n)}$,
  which is a gradient flow of the function
  $f:SO(n)\rightarrow\mathbb{R},\Theta\mapsto{2n-2\tr{\Theta}}$. Many of the
  presented results are not new, but scattered throughout literature. The
  motivation of this note is to summarize and extend the convergence results
  known from literature. Rather than giving an exhaustive list of references,
  the results are presented in a self-contained fashion.
\end{abstract}
\newcommand{\vectorize}{\mathop{\mathrm{vec}}}
\ForDiss{
In this section we discuss the properties of a function and a differential
equation on a smooth manifold.
}
\ForPaper{
In this note, we discuss the properties of a function and a differential
equation on a smooth manifold.
If we speak about a manifold~$\mathcal{M}$ of dimension $m$ we always mean a
smooth manifold in the
sense of \cite{Guillemin1974differentialTopology}, i.e. the subset of some
$\mathbb{R}^{k}$ with $k\geq{m}$ and $\mathcal{M}$ is locally diffeomorphic to
$\mathbb{R}^{m}$.
}
In the context of this note, we need the notions of \emph{measure zero} and \emph{dense}.
A set $A\subset\mathcal{M}$ of a manifold $\mathcal{M}$ is a set of
measure zero if there is a collection of smooth charts $\{U_l,\phi_l\}$ whose
domains cover $A$ and such that $\phi_l(A\cap{U_l})$ have
measure zero in $\mathbb{R}^{n}$, i.e. the $\phi_l(A\cap{U_l})$ can be
covered for any $\varepsilon$ by a countable collection of open balls whose
volumes sum up to less than $\varepsilon$, for details see e.g. \cite[Chapter
10]{Lee2006introductionToSmoothManifolds}.
To define \emph{dense}, we need the topological closure $\overline{A}$ of a set
$A\subset\mathcal{M}$, i.e. the intersection of all closed sets in $\mathcal{M}$
that contain $A$.
A dense subset of a smooth manifold $\mathcal{M}$ is a set
$A\subset\mathcal{M}$ such that the topological closure fulfills
$\overline{A}=\mathcal{M}$, see e.g. \cite[Appendix,
Topology]{Lee2006introductionToSmoothManifolds}.
$A$ is dense if and only if every nonempty open subset of $X$ has non-empty
intersection with $A$.
The complement $\mathbb{R}^{n}\setminus{A}$ of a 
set of measure zero $A\subset\mathbb{R}^{n}$ is dense in
$\mathbb{R}^{n}$, since if there is a point $x\in\mathbb{R}^{n}$ such that there
is an open $U\subset\mathbb{R}^{n}$ with $x\in{U}$ and
$U\cap{(\mathbb{R}^{n}\setminus{A})}=\emptyset$, then $A$
contains an open set and cannot have measure zero, see also
\cite[Chapter 2]{Milnor1997topologyFromTheDifferentiableViewpoint}.

Here, we consider a function and a differential equation on the set
of special orthogonal matrices
${SO(n)}=\{\Theta\in\mathbb{R}^{n\times{n}}\rvert{
\Theta^{-1}=\Theta^{\T},~\det(\Theta)=1}\}$.
$SO(n)$ is a smooth manifold of dimension $\tfrac{n(n-1)}{2}$ with the subspace
topology induced by $\mathbb{R}^{n\times{n}}$.
The tangent space $T_{\Theta}SO(n)$ at $\Theta$ is
given by
\begin{equation}
  \begin{aligned}
    T_{\Theta}SO(n)
    &=
    \{X\in\mathbb{R}^{n\times{n}}\rvert{
    X=\Omega\Theta,~
    \Omega\in\mathbb{R}^{n\times{n}},~
    \Omega=-\Omega^{\T}
    }\}
  \end{aligned}
  \text{.}
\end{equation}
The Riemannian metric
$g:T_{\Theta}SO(n)\times{T_{\Theta}SO(n)}\rightarrow\mathbb{R}$ induced by the
standard Euclidean metric on $SO(n)$ is given by
\begin{equation}
  \label{eqn:traceson:riemannianMetric}
  g(\Omega_1\Theta,\Omega_2\Theta)
  =
  \tr{(\Omega_1\Theta)^{\T}\Omega_2\Theta}
  =
  \tr{\Omega_1^{\T}\Omega_2}
  \text{.}
\end{equation}
In the following, we define the differential and the Hessian of a function
$f:SO(n)\rightarrow\mathbb{R}$ at a point $\Theta_0\in{SO(n)}$.
Let $\Gamma:(-\varepsilon,\varepsilon)\rightarrow{SO(n)}$ be a smooth curve with
$\dot{\Gamma}(t)=\Omega(t)\Gamma(t)$, $\Gamma(0)=\Theta_0$ and
$\dot{\Gamma}(0)=\Omega_0\Theta_0$ with $\Omega_0\in\mathbb{R}^{n\times{n}}$
and $\Omega_0=-\Omega_0^{\T}$.
The differential $\md{f}_{\Theta_0}:
T_{\Theta_0}SO(n)\rightarrow{T_{f(\Theta_0)}\mathbb{R}\cong\mathbb{R}}$ of a
function $f:SO(n)\rightarrow\mathbb{R}$ at a point $\Theta_0$ evaluated at
$\Omega_0\Theta_0\in{T_{\Theta_0}SO(n)}$ is defined by
\begin{equation}
  \md{f}_{\Theta_0}(\Omega_0\Theta_0)
  =
  \tfrac{\md}{\md{t}}\rvert_{t=0}(f\circ\Gamma)(t)
  \text{.}
  \label{}
\end{equation}
The \emph{critical points of $f$} are the points $\Theta_0$ where
$\md{f}_{\Theta_0}$ is not surjective.
Because of $\dim(T_{f(\Theta_0)}\mathbb{R})=1$, this means that
these are the points $\Theta_0$ where $\md{f}_{\Theta_0}=0$.
The \emph{gradient of $f$} is defined as the unique vector field $\grad{f}$ 
with
\begin{equation}
  df_{\Theta_0}(\Omega_0\Theta_0)
  =
  g(\grad{f}(\Theta_0),\Omega_0\Theta_0)
  \label{}
  \text{,}
\end{equation}
see e.g. \cite[Chapter 11]{Lee2006introductionToSmoothManifolds}.
The \emph{Hessian $H_f(\Theta_0)$ of $f$} at a critical point $\Theta_0$
evaluated at $(\Omega_0\Theta_0,\Omega_0\Theta_0)$ is defined by
\begin{equation}
  H_f(\Theta_0)(\Omega_0\Theta_0,\Omega_0\Theta_0)
  =
  \tfrac{\md^{2}}{\md{t}^{2}}\rvert_{t=0}(f\circ\Gamma)(t)
  \text{.}
  \label{}
\end{equation}
Since the Hessian at a critical point is bilinear and symmetric, we have
for $(\Omega_1+\Omega_2)\Theta_0\in{T_{\Theta_0}SO(n)}$ with
$\Omega_{1,2}=-\Omega_{1,2}^{\T}$ the equality
\begin{multline}
  \label{eqn:traceson:defn:hess:2}
  H_f(\Theta_0)((\Omega_1+\Omega_2)\Theta_0,(\Omega_1+\Omega_2)\Theta_0)
  \\
  =
  H_f(\Theta_0)(\Omega_1\Theta_0,\Omega_1\Theta_0)
  +
  2
  H_f(\Theta_0)(\Omega_1\Theta_0,\Omega_2\Theta_0)
  +
  H_f(\Theta_0)(\Omega_2\Theta_0,\Omega_2\Theta_0)
  \text{.}
\end{multline}
As a consequence, the value $H_f(\Theta_0)(\Omega_1\Theta_0,\Omega_2\Theta_0)$
can be computed utilizing the values
$H_f(\Theta_0)(\Omega_1\Theta_0,\Omega_1\Theta_0)$,
$H_f(\Theta_0)(\Omega_2\Theta_0,\Omega_2\Theta_0)$,
$H_f(\Theta_0)((\Omega_1+\Omega_2)\Theta_0,(\Omega_1+\Omega_2)\Theta_0)$
and \eqref{eqn:traceson:defn:hess:2}.
For details on the Hessian at a critical point, see
\cite[Appendix~C.5]{Helmke1996optimizationAndDynamicalSystems}.
\begin{lem}
  \label{lem:traceson:1}
  Consider the function
  $f:SO(n)\rightarrow\mathbb{R},\Theta\mapsto{n-\tr{\Theta}}$.
  \begin{enumerate}
    \item 
      \label{lem:traceson:1:1}
      The differential $\md{f}_{\Theta_0}$ of $f$ at $\Theta_0$ is given
      for any $\Omega_0\Theta_0\in{T_{\Theta_0}SO(n)}$ by
      \begin{equation*}
        \md{f}_{\Theta_0}(\Omega_0\Theta_0)
        =
        -\frac{1}{2}\tr{\Theta_0(\Theta_0-\Theta_0^{\T})\Omega_0\Theta_0}
        \label{}
      \end{equation*}
      and the critical points of $f$ are given by
      \begin{equation}
        \label{eqn:traceson:criticalSet}
        \mathcal{F}=\{\Theta_0\in{SO(n)}\rvert{\Theta_0^{\T}=\Theta_0}\}
        \text{.}
      \end{equation}
      Furthermore, the gradient~$\grad{f}(\Theta_0)$ at~$\Theta_0$
      is given by
      \begin{equation}
        \grad{f}({\Theta_0})
        =
        \frac{1}{2}(\Theta_0-\Theta_0^{\T})\Theta_0
        \text{.}
      \end{equation}
    \item
      \label{lem:traceson:1:2}
      The Hessian $H_{f}(\Theta_0)$ at a critical point $\Theta_0$ is given by
      \begin{equation*}
        H_f(\Theta_0)(\Omega_1\Theta_0,\Omega_2\Theta_0)
        =
        \frac{1}{2}\tr{
        \Omega_1^{\T}\Theta_0\Omega_2
        +
        \Omega_2^{\T}\Theta_0\Omega_1}
        \text{.}
        \label{}
      \end{equation*}
    \item
      \label{lem:traceson:1:3}
      The set of critical points $\mathcal{F}$ has the following properties:
      \begin{enumerate}
        \item 
          \label{lem:traceson:1:3:1}
          $\mathcal{F}=\cup_{k=0}^{\lfloor{\tfrac{n}{2}}\rfloor}\mathcal{F}_k$ where
          \begin{equation}
            \mathcal{F}_{k}
            =
            \{\Theta_0\in{SO(n)}\rvert{
            \Theta_0=\Theta_0^{\T},~
            \tr{\Theta_0}=n-4k
            }\}
            \text{.}
            \label{}
          \end{equation}
        \item 
          \label{lem:traceson:1:3:2}
          Each $\mathcal{F}_{k}$ is connected and isolated, i.e. there exists a
          neighborhood $U$ of each $\mathcal{F}_{k}$ such that
          $U\cap\mathcal{F}_{l}=\emptyset$ for all $l\neq{k}$.
        \item
          \label{lem:traceson:1:3:3}
          $\mathcal{F}_k$ are compact submanifolds of dimension $2k(n-2k)$
          and the tangent space at $\Theta_0\in\mathcal{F}_k$ is
          $T_{\Theta_0}\mathcal{F}_k=
          \{\Sigma\in\mathbb{R}^{n\times{n}}\rvert{\Sigma=\Sigma^{\T}}\}\cap{T_{\Theta_0}SO(n)}$.
        \item 
          \label{lem:traceson:1:3:4}
          For every $k\in\{0,\ldots,\lfloor{\tfrac{n}{2}}\rfloor\}$
          and every $\Theta_0\in\mathcal{F}_{k}$ we have
          \begin{multline*}
            \ker{H}_f(\Theta_0)
            =
            \{
            X\in{T_{\Theta_0}SO(n)}
            \rvert
            \\
            {
            \text{
            $H_f(\Theta_0)(X,Y)=0$
            for all
            $Y\in{T_{\Theta_0}SO(n)}$
            }
            }
            \}
            =
            T_{\Theta_0}\mathcal{F}_k
            \text{.}
          \end{multline*}
      \end{enumerate}
    \item
      \label{lem:traceson:1:4}
      $f$ has a unique minimum at $\Theta_0=I$, the other critical points are
      saddle points.
  \end{enumerate}
\end{lem}
\begin{cor}
  The differential equation
  \begin{equation}
    \dot{\Theta}
    =
    (\Theta^{\T}-\Theta)\Theta
    \label{eqn:traceson:ode}
  \end{equation}
  is the gradient flow of
  $f:SO(n)\rightarrow\mathbb{R},\Theta\mapsto{2n-2\tr{\Theta}}$ with respect to
  the Riemannian metric \eqref{eqn:traceson:riemannianMetric}.
\end{cor}
In the following we prove Lemma \ref{lem:traceson:1}.\\
\begin{mypf}
  \textbf{\ref{lem:traceson:1:1}}
  As stated above, $\Gamma:(-\varepsilon,\varepsilon)\rightarrow{SO(n)}$ is a
  differentiable curve with $\dot{\Gamma}(t)=\Omega(t)\Gamma(t)$,
  $\Gamma(0)=\Theta_0$ and $\dot{\Gamma}(0)=\Omega_0\Theta_0$ with
  $\Omega_0\in\mathbb{R}^{n\times{n}}$ and $\Omega_0=-\Omega_0^{\T}$.
  Then
  \begin{equation}
    \begin{aligned}
      \md{f}_{\Theta_{0}}(\Omega_0\Theta_{0})
      &=
      \frac{\md}{\md{t}}\rvert_{t=0}
      \left(
      n-\tr{\Gamma(t)}
      \right)
      =
      -\tr{\Omega_0\Theta_{0}}
      \\
      &=
      -\frac{1}{2}\tr{\Theta_{0}^{\T}\Omega_0^{\T}+\Omega_0\Theta_{0}}
      \\
      &=
      -\frac{1}{2}\tr{(\Theta_{0}-\Theta_{0}^{\T})\Omega_0}
      \\
      &=
      -\frac{1}{2}
      \tr{\Theta_{0}^{\T}(\Theta_{0}-\Theta_{0}^{\T})\Omega_0\Theta_{0}}
      \text{.}
    \end{aligned}
  \end{equation}
  Therefore, the critical points of~$f$ are given by
  \begin{equation}
    \{\Theta_0\in{SO(n)}\rvert{
    \Theta_0=\Theta_0^{\T}}\}
    \text{.}
  \end{equation}
  With the definition of the Riemannian metric by \eqref{eqn:traceson:riemannianMetric},
  the gradient $\grad{f}_{\Theta}$ at~$\Theta_0$ is given by
  \begin{equation}
    \grad{f}(\Theta_0)
    =
    \frac{1}{2}(\Theta_{0}-\Theta_{0}^{\T})\Theta_{0}^{\T}
    \text{.}
  \end{equation}
  \textbf{\ref{lem:traceson:1:2}}
  Let $\Theta_0$ denote a critical point of $f$.
  As stated above, $\Gamma:(-\varepsilon,\varepsilon)\rightarrow{SO(n)}$ is a
  differentiable curve with $\dot{\Gamma}(t)=\Omega(t)\Gamma(t)$,
  $\Gamma(0)=\Theta_0$ and $\dot{\Gamma}(0)=\Omega_0\Theta_0$ with
  $\Omega_0\in\mathbb{R}^{n\times{n}}$ and $\Omega_0=-\Omega_0^{\T}$.
  Then
  \begin{equation}
    \begin{aligned}
      H_f(\Theta_0)(\Omega_0\Theta_0,\Omega_0\Theta_0)
      &=
      \frac{\md^{2}}{\md{t^{2}}}\rvert_{t=0}(f(\Gamma(t))
      =
      -
      \tr{\ddot{\Gamma}(t)}\rvert_{t=0}
      \\
      &=
      -\tr{\dot{\Omega}(t)\Gamma(t)+\Omega(t)\dot{\Gamma}(t)}
      \rvert_{t=0}
      \\
      &=
      -\tr{\Omega_0^{2}\Theta_0}
      =
      \tr{\Omega_0^{\T}\Theta_0\Omega}
      =
      \tr{\Theta_0^{\T}\Omega_0^{\T}\Theta_0\Omega_0\Theta_0}
      \text{,}
    \end{aligned}
    \label{eqn:traceson:computationHessian}
  \end{equation}
  where we utilized
  that~$\tr{\dot{\Omega}(0)\Theta_0}=0$ since~$\dot{\Omega}(t)$ is skew
  symmetric for all $t$ and $\Theta_0=\Theta_0^{\T}$ since $\Theta_0$ is a
  critical point.
  Utilizing \eqref{eqn:traceson:defn:hess:2} we get
  \begin{equation}
    \begin{aligned}
      H_f(\Theta_{0})(\Omega_1\Theta_{0},\Omega_2\Theta_{0})
      &=
      \frac{1}{2}\biggl(
      H_f(\Theta_{0})((\Omega_1+\Omega_2)\Theta_{0},(\Omega_1+\Omega_2)\Theta_{0})
      \\
      &
      ~
      ~
      ~
      ~
      -
      H_f(\Theta_{0})(\Omega_1\Theta_{0},\Omega_1\Theta_{0})
      -
      H_f(\Theta_{0})(\Omega_2\Theta_{0},\Omega_2\Theta_{0})
      \biggr)
      \\
      &=
      \frac{1}{2}\tr{
      \!\!
      (\Omega_1+\Omega_2)^{\T}\Theta_{0}(\Omega_1+\Omega_2)
      \!-\!
      \Omega_1^{\T}\Theta_{0}\Omega_1
      \!-\!
      \Omega_2^{\T}\Theta_{0}\Omega_2
      \!\!
      }
      \\
      &=
      \frac{1}{2}\tr{
      \Omega_1^{\T}\Theta_{0}\Omega_2
      +
      \Omega_2^{\T}\Theta_{0}\Omega_1
      }
      \text{.}
    \end{aligned}
    \label{}
  \end{equation}
  ~\\
  \textbf{\ref{lem:traceson:1:3:1}}
  Since~$\Theta_0$ is symmetric, $\Theta_0$ is orthonormally diagonalizable,
  i.e. $\Theta_0=\Pi^{\T}D\Pi$ for some diagonal~$D$ and orthonormal~$\Pi$ where
  the columns of~$\Pi$ are eigenvectors of~$\Theta_0$.
  Since~$\Theta_0^{\T}\Theta_0=I$, we get~$D^2=I$ and consequently the
  eigenvalues are~$\pm{1}$.  Since~$\Theta_0\neq{I}$ and~$\det(\Theta_0)=1$, we
  always have an even number of negative eigenvalues.
  A similarity transformation leaves the trace
  invariant, hence a critical point $\Theta_{0}$ fulfills
  \begin{equation}
    \tr{\Theta_{0}}=n-4k
    \text{,}
  \end{equation}
  where $k\in\{0,\ldots,\lfloor{\frac{n}{2}}\rfloor\}$ is the number of
  eigenvalue pairs which are $-1$.
  ~\\
  \textbf{\ref{lem:traceson:1:3:2}}
  We start by showing that each $\mathcal{F}_k$ is path
  connected and thus connected.
  Let $k\in\{0,\ldots,\lfloor{\frac{n}{2}}\rfloor\}$ be arbitrary but fixed
  and let $\Theta_1,\Theta_2\in\mathcal{F}_k$.
  Then there are orthogonal $\Pi_1,\Pi_2$ such that $\Theta_1=\Pi_1^{\T}D\Pi_1$
  and $\Theta_2=\Pi_2^{\T}D\Pi_2$. Furthermore, there are real
  skew-symmetric matrices $\Omega_1,\Omega_2$ such that $\Pi_1=\exp(\Omega_1)$
  and $\Pi_2=\exp(\Omega_2)$ with $\exp$ denoting the matrix exponential.
  Then $\alpha:[0,1]\rightarrow{SO(n)}$ defined by
  \begin{equation}
    t
    \mapsto
    \exp\left(\Omega_1^{\T}(1-t)\right)
    \exp\left(\Omega_2^{\T}t\right)
    D
    \exp\left(\Omega_2t\right)
    \exp\left(\Omega_1(1-t)\right)
  \end{equation}
  is a smooth curve in $\mathcal{F}_k$ which connects $\Theta_1$ and $\Theta_2$.
  Since $\Theta_1,\Theta_2\in\mathcal{F}_k$ were arbitrary, this implies the
  path-connectedness of $\mathcal{F}_k$.
  To show that $\mathcal{F}_k$ is isolated, we utilize that
  $n-4l=f\rvert_{\mathcal{F}_l}\neq{f\rvert_{\mathcal{F}_k}}=n-4k$ for
  $l\neq{k}$.
  Then there is a $\varepsilon(l)$ with
  $(n-4l-\varepsilon(l),n-4l+\varepsilon(l))
  \cap(n-4k-\varepsilon(l),n-4k+\varepsilon(l))
  =\emptyset$ for $k\neq{l}$.  As a consequence, the intersection of the
  preimage of these sets under $f$ is empty.
  Since $f$ is continuous and both, $(n-4l-\varepsilon(l),n-4l+\varepsilon(l))$
  and $(n-4k-\varepsilon(l),n-4k+\varepsilon(l))$ are open, their preimages are
  open and contain $\mathcal{F}_l$ and $\mathcal{F}_k$ respectively.
  With $U_k(l)=f^{-1}\bigl((n-4l-\varepsilon(l),n-4l+\varepsilon(l))\bigr)$
  we thus have $U_k(l)\cap{\mathcal{F}_l}=\emptyset$.
  Since this is possible for every
  $l\in\{0,\ldots,\lfloor{\tfrac{n}{2}}\rfloor\}$ and since a finite
  intersection of open sets is an open set, we find an open neighborhood $U$ of
  $\mathcal{F}_k$ such that $U\cap{\mathcal{F}_l}=\emptyset$ for all
  $l\in\{0,\ldots,\lfloor{\tfrac{n}{2}}\rfloor\}$.
  ~\\
  \textbf{\ref{lem:traceson:1:3:3}}
  The property that the $\mathcal{F}_k$ are submanifolds is given in
  \cite{Frankel1965criticalSubmanifoldsOfTheClassicalGroupsAndStiefelManifolds}.
  The tangent space follows from \eqref{eqn:traceson:criticalSet}.
  ~\\
  \textbf{\ref{lem:traceson:1:3:4}}
  Let $k\in\{0,\ldots,\lfloor{\frac{n}{2}}\rfloor\}$ be arbitrary but fixed and
  $\Theta_{0}\in\mathcal{F}_k$.
  Since
  $\ker{H}_f(\Theta_{0})\supset{T_{\Theta_{0}}\mathcal{F}_k}$ is always true,
  we have to check
  $\ker{H}_f(\Theta_{0})\subset{T_{\Theta_{0}}\mathcal{F}_k}$.  Since every
  critical point is symmetric, there is an orthogonal $\Pi$ such that
  $\Pi^{\T}\Theta_{0}\Pi=D$ where $D$ is a diagonal matrix with non-zero
  diagonal elements.  Furthermore, we know that
  \begin{equation*}
    \begin{aligned}
      H_f(\Theta_{0})(\Omega_1\Theta_{0},\Omega_2\Theta_{0})
      &=
      -\frac{1}{2}\tr{
      \Omega_1^{\T}\Theta_{0}\Omega_2
      +
      \Omega_2^{\T}\Theta_{0}\Omega_1
      }
      \\
      &=
      -\frac{1}{2}\tr{
      \Theta_{0}^{\T}\Omega_1^{\T}\Theta_{0}\Omega_2\Theta_{0}
      +
      \Theta_{0}^{\T}\Omega_2^{\T}\Theta_{0}\Omega_1\Theta_{0}
      }
      \\
      &=
      -\frac{1}{2}\tr{
      (\Omega_1\Pi^{\T}D\Pi)^{\T}\Pi^{\T}D\Pi(\Omega_2\Pi^{\T}D\Pi)
      }
      \\
      &
      ~
      ~
      ~
      -\frac{1}{2}\tr{
      (\Omega_2\Pi^{\T}D\Pi)^{\T}\Pi^{\T}D\Pi(\Omega_1\Pi^{\T}D\Pi)
      }
      \\
      &=
      -\frac{1}{2}\tr{
      (\Pi\Omega_1\Pi^{\T}D)^{\T}D(\Pi\Omega_2\Pi^{\T}D)
      }
      \\
      &
      ~
      ~
      ~
      -\frac{1}{2}\tr{
      (\Pi\Omega_2\Pi^{\T}D)^{\T}D(\Pi\Omega_1\Pi^{\T}D)
      }
      \\
      &=
      -\frac{1}{2}\tr{
      (\tilde{\Omega}_1D)^{\T}D(\tilde{\Omega}_2{D})
      +
      (\tilde{\Omega}_2D)^{\T}D(\tilde{\Omega}_1{D})
      }
      \\
      &=
      H_{f}(D)(\tilde{\Omega}_1D,\tilde{\Omega}_2D)
    \end{aligned}
  \end{equation*}
  where $\tilde{\Omega}_1=\Pi\Omega_1\Pi^{\T}$ and
  $\tilde{\Omega}_2=\Pi\Omega_2\Pi^{\T}$.
  As a consequence
  \begin{equation*}
    \begin{aligned}
      \ker{H}_f(\Theta_{0})
      &=
      \{
      X\in{T_{\Theta_{0}}SO(n)}
      \rvert{
      \text{
      $H_f(\Theta_{0})(X,Y)=0$
      for all
      $Y\in{T_{\Theta_{0}}SO(n)}$
      }
      }
      \}
      \\
      &=
      \{
      X\in{T_{D}SO(n)}
      \rvert{
      \text{
      $
      \tr{
      X^{\T}D\tilde{\Omega}_2D
      -
      D\tilde{\Omega}_2DX
      }
      =
      0
      $
      for all
      $\tilde{\Omega}_2{D}\in{T_{D}SO(n)}$
      }
      }
      \}
      \\
      &=
      \{
      X\in{T_{D}SO(n)}
      \rvert{
      \text{
      $
      \tr{
      \tilde{\Omega}_2D(X^{\T}-X)D
      }
      =
      0
      $
      for all
      $\tilde{\Omega}_2{D}\in{T_{D}SO(n)}$
      }
      }
      \}
      \text{.}
    \end{aligned}
  \end{equation*}
  Observe that $\tilde{\Omega}_2$ is skew symmetric.
  Furthermore, $X^{\T}-X$ is skew symmetric and since $D$ is diagonal with
  non-zero entries, $D(X^{\T}-X)D$ is skew symmetric.
  Because the equation $\tr{D\tilde{\Omega}_2D(X^{\T}-X)}=0$ has to hold for all
  skew-symmetric $\tilde{\Omega}_2$, we obtain $D(X^{\T}-X)D=0$.
  With the non-singular $D$ this implies $X^{\T}-X=0$ which is equivalent to
  $X=X^{\T}$.
  In \ref{lem:traceson:1:3:2} we showed $T_{\Theta_0}\mathcal{F}_k$ is
  $\{\Sigma\in\mathbb{R}^{n\times{n}}\rvert{\Sigma=\Sigma^{\T}}\}\cap{T_{\Theta_0}SO(n)}$, thus
  the previous calculation shows
  $\ker{H}_{f}(\Theta_0)\subset{T_{\Theta_0}\mathcal{F}_k}$.
  ~\\
  \textbf{\ref{lem:traceson:1:4}}
  Since $\Theta_{0}=\Theta_{0}^{\T}$, $\Theta_{0}$ is orthogonally diagonalizable,
  i.e.~$\Theta_{0}=\Pi^{\T}D\Pi$ for some diagonal~$D$ and orthogonal~$\Pi$.
  Therefore
  \begin{equation}
    \begin{aligned}
      H_{f}((\Omega_{0}\Theta_{0}),(\Omega_{0}\Theta_{0}))
      &=
      \tr{\Omega_{0}^{\T}\Theta_{0}\Omega_{0}}
      =
      -\tr{\tilde{\Omega}_{0}^{2}D}
      \text{,}
    \end{aligned}
  \end{equation}
  where $\tilde{\Omega_{0}}=\Pi\Omega_{0}\Pi^{\T}$ is skew symmetric.
  Consequently, $\tr{\Omega_{0}^{\T}\Theta_{0}\Omega_{0}}$ is definite for all skew
  symmetric $\Omega_{0}$ at a critical point $\Theta_{0}$ if and only
  if $\tr{\tilde{\Omega_{0}}^{2}D}$ is definite for all skew
  symmetric $\tilde{\Omega_{0}}$ where $D=\Pi\Theta_{0}\Pi^{\T}$ is diagonal. Thus,
  we have to consider $H_{f}$ only for diagonal $D$, i.e.
  \begin{equation}
    \begin{aligned}
      H_{f}((\Omega_{0}\Theta_{0}),(\Omega_{0}\Theta_{0}))
      &=
      -\tr{\Omega_{0}^{2}{D}}
      =
      -\sum_{1\leq{k}\leq{n}}{\left(\Omega_{0}^{2}D\right)_{kk}}
      \\
      &=
      -\sum_{1\leq{k}\leq{n}}{
      (\Omega_{0}^{2})_{kk}D_{kk}
      }
      =
      -\sum_{1\leq{k}\leq{n}}{
      \sum_{1\leq{l}\leq{n}}{
      (\Omega_{0})_{kl}(\Omega_{0})_{lk}D_{kk}
      }
      }
      \\
      &=
      \sum_{1\leq{k,l}\leq{n}}{
      ((\Omega_{0})_{kl})^{2}D_{kk}
      }
      =
      \sum_{
      \substack{
      1\leq{k,l}\leq{n}\\
      k\neq{l}
      }
      }
      {
      ((\Omega_{0})_{kl})^{2}D_{kk}
      }
      \\
      &=
      \sum_{1\leq{k}<{l}\leq{n}}{
      (D_{kk}+D_{ll})((\Omega_{0})_{kl})^2
      }
      \text{.}
    \end{aligned}
    \label{eqn:traceson:1}
  \end{equation}
  The critical points $\Theta_{0}$ are such that $\Theta_{0}\in{SO(n)}$ are
  symmetric, therefore all eigenvalues of $\Theta_{0}$ are real.
  Observe now that $\Theta_{0}=\Pi^{\T}D\Pi$ with orthogonal $\Pi$
  implies $D^{2}=I$.
  Hence, the eigenvalues are $D_{kk}\in\{-1,1\}$ and
  since $\Theta_{0}\in{SO(n)}$ the number of $-1$-eigenvalues is even.
  Consequently we have to determine the definiteness
  of $H_{f}$ by considering \eqref{eqn:traceson:1} for all
  diagonal matrices $D$ with $\pm{1}$ on the diagonals where the number of $-1$
  entries is zero or even.
  ~\\
  Suppose first, that all $D_{kk}$ are equal to $1$, i.e. $D=I$.
  The associated $\Theta_{0}$ is $\Theta_{0}=\Pi^{\T}D\Pi=\Pi^{\T}\Pi=I$.
  $\eqref{eqn:traceson:1}$ then implies
  $H_{f}(\Omega_{0}\Theta_{0},\Omega_{0}\Theta_{0})
  =
  \sum_{1\leq{k}<{l}\leq{n}}{2(\Omega_{0})_{kl}^2}$, i.e. $H_{f}>0$ for all
  skew-symmetric $\Omega_{0}$. 
  Thus $H_{f}$ is positive definite if $\Theta_{0}=I$.
  Suppose now there is an even number of eigenvalues $D_{kk}$ equal to $-1$.
  Then, there are
  indices $l,k$ and $l\neq{k}$ such that $D_{kk}=-1$ and $D_{ll}=-1$, and
  therefore there are skew symmetric $\Omega_{0}$ such
  that $H_{f}(\Omega_{0}\Theta_{0},\Omega_{0}\Theta_{0})<0$.
  As consequence, $H_{f}$ is indefinite at a critical
  point $\Theta_{0}$ where $\Theta_{0}$ has an even number of
  negative eigenvalues $D_{kk}=-1$.
  Therefore, $\Theta_{0}=I$ is the only local (global) minimum
  of $f$. All other critical points are saddle points.
\end{mypf}
\begin{defn}\cite[on p.21]{Helmke1996optimizationAndDynamicalSystems}
  \label{defn:traceson:3}
  Let $\mathcal{M}$ be a smooth Riemannian manifold and
  $f:\mathcal{M}\rightarrow\mathbb{R}$ be a smooth function. Denote the set of
  critical points of $f$ by $C(f)$. $f$ is called \emph{Morse-Bott function}
  provided the following conditions are satisfied:
  \begin{enumerate}
    \item
      \label{defn:traceson:3:1}
      $f$ has compact sublevel sets.
    \item
      \label{defn:traceson:3:2}
      $C(f)=\cup_{j=1}^{k}\mathcal{N}_j$ where $\mathcal{N}_j$ are disjoint,
      closed and connected submanifolds of $\mathcal{M}$ and $f$ is constant on
      $\mathcal{N}_j$ for $j=1,\ldots,k$.
    \item
      \label{defn:traceson:3:3}
      $\ker{H_{f}(x)}=T_x\mathcal{N}_k$ for all $x\in\mathcal{N}_j$ and
      all $j=1,\ldots,k$.
  \end{enumerate}
\end{defn}
\begin{lem}
  \label{lem:traceson:4}
  $f:SO(n)\rightarrow\mathbb{R},\Theta\mapsto{n-\tr{\Theta}}$
  is a Morse-Bott function.
\end{lem}
\begin{mypf}
  We show only Definition~\ref{defn:traceson:3}\ref{defn:traceson:3:1} since \ref{defn:traceson:3:2} and
  \ref{defn:traceson:3:3} were shown in Lemma \ref{lem:traceson:1}.
  $SO(n)$ is compact, hence $f$ attains its minimal and its maximal value on
  $SO(n)$. The minimal value of $f$ is zero, the maximal value is $2n-2$ for $n$
  odd and $2n$ for $n$ even.
  If $n$ is odd we thus have
  \begin{equation*}
    L_{c}
    =\{\Theta\in{SO(n)}\rvert{
    f(\Theta)\leq{c}
    }\}
    =
    \begin{cases}
      f^{-1}([0,c])
      &
      \text{for $c\leq{2n-2}$}
      \\
      f^{-1}([0,2n-2])
      &
      \text{for $c>{2n-2}$}
      \text{.}
    \end{cases}
  \end{equation*}
  If $n$ is even we have
  \begin{equation*}
    L_{c}
    =\{\Theta\in{SO(n)}\rvert{
    f(\Theta)\leq{c}
    }\}
    =
    \begin{cases}
      f^{-1}([0,c])
      &
      \text{for $c\leq{2n}$}
      \\
      f^{-1}([0,2n])
      &
      \text{for $c>{2n}$}
      \text{.}
    \end{cases}
  \end{equation*}
  Since $f$ is continuous, the preimage of a closed set is a closed set and
  since $SO(n)$ is bounded, its subsets are bounded as well. Since
  $SO(n)\subset\mathbb{R}^{n\times{n}}$, the boundedness and closedness of the
  sublevel sets implies their compactness.
\end{mypf}
The convergence properties of the gradient flow associated with
$\Theta\mapsto{n-\tr{\Theta}}$ are thus determined by the following proposition.
\begin{prop}\cite[Proposition 3.9]{Helmke1996optimizationAndDynamicalSystems}
  \label{prop:traceson:5}
  Let $f:\mathcal{M}\rightarrow\mathbb{R}$ be a Morse-Bott function on a
  Riemannian manifold $\mathcal{M}$. The $\omega$-limit set $\omega(x)$ of
  $x\in\mathcal{M}$ with respect to the gradient flow of $f$ is a single
  critical point of $f$.
  Every solution of the gradient flow converges to an equilibrium point.
\end{prop}

To give a more detailed specification the convergence behavior of the gradient
flow of $\Theta\mapsto{n-\tr{\Theta}}$ we need the following result.

\begin{lem}
  \label{lem:traceson:6}
  Let $\mathcal{M}$ be a smooth and compact Riemannian manifold of dimension
  $m$, $f:\mathcal{M}\rightarrow\mathbb{R}$ be a Morse-Bott function and denote
  the set of critical points of $f$ by $C(f)$.
  Let $\mathcal{N}$ be a fixed connected component of $C(f)$
  of dimension $n$.
  If at least one of the $m-n$ eigenvalues with nonzero real part of the
  linearization of $\grad{f}$ at some $x\in\mathcal{N}$ has a real part greater than
  zero, then the set ${A}$ of initial conditions $x_0\in\mathcal{M}$
  for which the solutions $t\mapsto{\phi(t,x_0)}$ of the gradient flow
  $\dot{x}=-\grad{f}(x)$ converge towards $\mathcal{N}$, i.e.
  \begin{equation}
    {A}
    =
    \{x_0\in\mathcal{M}
    \rvert{\lim_{t\rightarrow\infty}\phi(t,x_0)\in\mathcal{N}}
    \}
    \text{,}
    \label{eqn:traceson:inset}
  \end{equation}
  has measure zero. Furthermore $\mathcal{M}\setminus{{A}}$ is dense in
  $\mathcal{M}$, i.e.  $\overline{M\setminus{{A}}}=\mathcal{M}$.
\end{lem}
\begin{mypf}
  The goal of the proof is to show that ${A}$ has measure zero and that
  $\mathcal{M}\setminus{A}$ is dense. We show this in the following way.
  First, we consider the set of points lying in a suitable neighborhood of
  $\mathcal{N}$ and which contains the orbits of the
  solutions of the gradient flow
  $\dot{x}=-\grad{f}(x)$ which eventually converge towards $\mathcal{N}$.
  We utilize a result from
  \cite{Aulbach1984continuousAndDiscreteDynamcisNearManifoldsOfEquilibria} to
  conclude that this set has measure zero and $\mathcal{M}$ without this set is
  dense. Then, we utilize this set to derive the same result for ${A}$
  utilizing the properties of the flow of the gradient vector field on
  $\mathcal{M}$.

  In the following, we apply \cite[Proposition
  4.1]{Aulbach1984continuousAndDiscreteDynamcisNearManifoldsOfEquilibria}.
  This proposition concerns the case of a a three times continuously
  differentiable vector field $v:\mathbb{R}^{l}\rightarrow\mathbb{R}^{l}$
  together with 
  submanifold of equilibria $\overline{\mathcal{N}}$ in $\mathbb{R}^{l}$
  under the assumption that
  $\overline{\mathcal{N}}$ is normally hyperbolic with respect to $v$.
  Normal hyperbolicity of $\overline{\mathcal{N}}$ means that the linearization
  of the vector field $v$ at $x\in\overline{\mathcal{N}}$ has
  $n-\dim{\overline{\mathcal{N}}}$ eigenvalues with real parts different from
  zero.
  Under these assumptions, there exists a neighborhood
  $\overline{\mathcal{U}}$ of $\overline{\mathcal{N}}$ such that any solution
  $t\mapsto\phi(t,x_0)$ of $\dot{x}=v(x)$ with initial condition $x_0$ and with
  a forward
  orbit $\phi([0,\infty);x_0)$ in $\overline{\mathcal{U}}$ lies on the stable of
  manifold $W^{s}_{\text{loc}}(p)$ of a point $p\in\overline{\mathcal{N}}$.
  $W^{s}_{\text{loc}}(p)$ is defined by 
  \begin{equation}
    W^{s}_{\text{loc}}(p)
    =
    \{x\in\mathcal{U}\rvert{\lim_{t\rightarrow\infty}\phi(t,x)=p}\}
    \text{.}
    \label{}
  \end{equation}
  We can always embed $\mathcal{M}$ into $\mathbb{R}^{l}$ for $l$ large enough,
  see e.g. \cite[Chapter 10]{Lee2006introductionToSmoothManifolds}, therefore we
  can utilize
  \cite{Aulbach1984continuousAndDiscreteDynamcisNearManifoldsOfEquilibria} also
  for our case of a vector field on a manifold $\mathcal{M}$.

  Since $\mathcal{M}$ is compact and $f$ is smooth, we have a global flow
  $\phi:\mathbb{R}\times\mathcal{M}\rightarrow\mathcal{M}$, which means that
  $t\mapsto{\phi(t,x_0)}$ is a solution of the gradient flow
  $\dot{x}=-\grad{f}(x)$ defined for all $t\in\mathbb{R}$ and with
  $\phi(0,x_0)=x_0$, see e.g.  \cite[Chapter
  17]{Lee2006introductionToSmoothManifolds}.
  Furthermore $\phi(t,\cdot):\mathcal{M}\rightarrow\mathcal{M}$ is a
  diffeomorphism for every $t\in\mathbb{R}$.
  Since $f$ is a Morse-Bott function, $\mathcal{N}$ is normally hyperbolic,
  i.e. the linearization of the gradient flow at any $x\in\mathcal{N}$ 
  has exactly $m-n$ eigenvalues with real parts different from zero,
  see \cite[p. 183, Morse-Bott
  functions]{McDuff1998introductionToSymplecticTopology}.
  According to 
  \cite[Proposition
  4.1]{Aulbach1984continuousAndDiscreteDynamcisNearManifoldsOfEquilibria}, we
  have a neighborhood $\mathcal{U}$ of $\mathcal{N}$ such that for every
  solution $\phi(\cdot,x)$ with $x\in\mathcal{N}$ and a forward orbit
  $\phi([0,\infty],x)$ in $\mathcal{U}$, the solution has to lie in one
  $W^{s}_{\text{loc}}(p)$ with $p\in\mathcal{N}$.
  We know from \cite[Proposition
  3.2]{Austin1995morseBottTheoryAndEquivariantCohomology},
  that if we choose $\mathcal{U}$ small
  enough, then the local stable manifold
  $W^{s}_{\text{loc}}(\mathcal{N})$ of $\mathcal{N}$ given by
  \begin{equation}
    W^{s}_{\text{loc}}(\mathcal{N})
    =
    \cup_{p\in\mathcal{N}}W^{s}_{\text{loc}}(p)
    \label{}
  \end{equation}
  is a smooth submanifold of dimension $m+k$ where $k<m-n$ is the number of
  eigenvalues with real part smaller than zero.
  Since the stable manifold $W^{s}_{\text{loc}}(\mathcal{N})$ is a submanifold
  of $\mathcal{M}$ with smaller dimension than $\mathcal{M}$, the stable
  manifold has measure zero and $\mathcal{M}\setminus{W^{s}_{\text{loc}}(\mathcal{N})}$
  is dense in $\mathcal{M}$, see \cite[Theorem
  10.5]{Lee2006introductionToSmoothManifolds}.

  Let ${A}$ be defined by \eqref{eqn:traceson:inset}.
  Define ${A}_1$ by
  \begin{equation}
    {A}_1
    =
    \{x\in{A}\rvert{
    \forall{t\geq{1}}:
    \phi(t,x)\in\mathcal{U}
    }\}
    \label{}
  \end{equation}
  and let ${A}_k$ for $k\geq{2}$ be defined by
  \begin{equation}
    {A}_k
    =
    \{
    x
    \in
    {A}\setminus\bigl({A}_1\cup\ldots\cup{A}_{k-1}\bigr)
    \rvert{
    \forall{t\geq{k}}:
    \phi(t,x)\in\mathcal{U}
    }\}
    \text{.}
    \label{eqn:traceson:defn:mathcalA:2}
  \end{equation}
  If $x\in{A}$, then there is an integer $k\in\mathbb{N}$ such that
  $x\in{A}_k$. As a consequence
  \begin{equation}
    {A}
    =
    \cup_{k\in\mathbb{N}}
    {{A}_k}
    \text{.}
    \label{eqn:traceson:prop:mathcalA}
  \end{equation}
  Because of \eqref{eqn:traceson:defn:mathcalA:2},
  $\phi(k,{A}_k)\subset\mathcal{U}$ for every ${A}_k$.
  Moreover, 
  \cite[Proposition
  4.1]{Aulbach1984continuousAndDiscreteDynamcisNearManifoldsOfEquilibria}
  implies that
  \begin{equation}
    \phi(k,{A}_k)
    \subset
    {W^{s}_{\text{loc}}(\mathcal{N})}
    \text{.}
    \label{eqn:traceson:prop:mathcalAk}
  \end{equation}
  As subset of a set of measure zero, $\phi(k,{A}_k)$ has measure zero,
  see e.g.  \cite[Lemma A.60(b)]{Lee2006introductionToSmoothManifolds}. Since
  $\phi(k,\cdot):\mathcal{M}\rightarrow\mathcal{M}$ is a diffeomorphism, this
  means that ${A}_k$ also has measure zero, see e.g. \cite[Lemma
  10.1]{Lee2006introductionToSmoothManifolds}. According to
  \eqref{eqn:traceson:prop:mathcalA}, ${A}$ is a countable union of the
  ${A}_k$, i.e. ${A}$ is a countable union of sets of measure
  zero. Therefore, ${A}$ has measure zero and as a consequence
  $\mathcal{M}\setminus{A}$ is dense.

\end{mypf}

To finally derive the global stability properties of the identity matrix
$I\in\mathbb{R}^{n\times{n}}$ for the gradient flow of
$\Theta\mapsto{n-\tr{\Theta}}$ and thus also for the differential equation
\eqref{eqn:traceson:ode}, we linearize the gradient flow around the
equilibria.

\begin{lem}
  \label{lem:traceson:7}
  The convergence properties of the gradient flow of
  $f:SO(n)\rightarrow\mathbb{R},\Theta\mapsto{n-\tr{\Theta}}$ are the following:
  \begin{enumerate}
    \item
      \label{lem:traceson:7:1}
      The $\omega$-limit set of any solution is contained in the set of
      equilibria given by \eqref{eqn:traceson:criticalSet}, i.e.
      \begin{equation*}
        \mathcal{F}=\{\Theta_0\in{SO(n)}\rvert{\Theta_0^{\T}=\Theta_0}\}
        \text{.}
      \end{equation*}
    \item
      \label{lem:traceson:7:2}
      The equilibrium~$I$ is locally exponentially stable
      and all other equilibria are unstable.
    \item
      \label{lem:traceson:7:3}
      The set of initial conditions for which the solutions of the gradient flow
      $\dot{x}=-\grad{f(x)}$ of $f$ converge towards $I$ is dense in $SO(n)$
      and the set of initial conditions for which the solutions of the gradient
      flow converge to the other equilibria has measure zero.
  \end{enumerate}
\end{lem}
\begin{cor}
  \label{cor:traceson:8}
  The identity matrix $I$ is an almost globally asymptotically stable
  equilibrium for the differential equation
  \begin{equation}
    \dot{\Theta}
    =
    (\Theta^{\T}-\Theta)\Theta
    \text{.}
  \end{equation}
\end{cor}
In the following we prove Lemma \ref{lem:traceson:7}.\\
\begin{mypf}
  \ref{lem:traceson:7:1} is a consequence of Lemma \ref{lem:traceson:4} and
  Proposition \ref{prop:traceson:5}.\\
  \ref{lem:traceson:7:2}
  To prove the property \ref{lem:traceson:7:2} we
  linearize the gradient flow with the vector field $\grad{f}$ defined by
  $\Theta\mapsto\tfrac{1}{2}(\Theta^{\T}-\Theta)\Theta$ around the equilibria.
  To do this directly we compute
  $\md{\grad{f}}_{\Theta_0}(X)
  =\tfrac{\md}{\md{t}}\rvert_{t=0}((\grad{f})\circ{\Gamma})(t)$
  where $\Theta_0$ is an equilibrium and 
  $\Gamma:(-\varepsilon,\varepsilon)\rightarrow{SO(n)}$ is smooth with
  $\Gamma(0)=\Theta_0$, $\dot{\Gamma}(0)=X$ and $X\in{T_{\Theta_0}}SO(n)$.
  This yields
  \begin{equation}
    \begin{aligned}
      {\md\grad{f}}_{\Theta_0}(X)
      &=
      \frac{1}{2}
      \frac{\md}{\md{t}}\rvert_{t=0}
      (
      \Gamma^{\T}(t)-\Gamma(t)
      )\Gamma(t)
      \\
      &=
      \frac{1}{2}\left(
      \dot{\Gamma}^{\T}(t)\Gamma(t)
      +
      \Gamma^{\T}(t)\dot{\Gamma}(t)
      -
      \dot{\Gamma}(t)\Gamma(t)
      -
      \Gamma(t)\dot{\Gamma}(t)
      \right)\rvert_{t=0}
      \\
      &=
      \frac{1}{2}(
      X^{\T}\Theta_0+\Theta_0^{\T}X-X\Theta_0-\Theta_0X
      )
      \\
      &=
      -\frac{1}{2}\left(
      \Theta_0X
      +
      X\Theta_0
      \right)
      \text{,}
    \end{aligned}
  \end{equation}
  since $\Theta_0=\Theta_0^{\T}$ and $X=\Omega_0\Theta_0$ for a
  $\Omega_0\in\mathbb{R}^{n\times{n}}$ with $\Omega_0=-\Omega_0^{\T}$.
  Consequently the linearization of the gradient flow at an equilibrium
  $\Theta_0$ is given by
  \begin{equation}
    \begin{aligned}
      \dot{X}
      &=
      -\frac{1}{2}
      \left(
      \Theta_0^{\T}X
      +
      X{\Theta_0}
      \right)
      \text{,}
    \end{aligned}
    \label{eqn:traceson:linearizationGradient}
  \end{equation}
  where~$X\in{T_{\Theta_0}SO(n)}$. Note that due to the simple nature of the
  Riemannian metric \eqref{eqn:traceson:riemannianMetric} and the connection of
  the linearization of a gradient flow to the Hessian, we could have obtained
  \eqref{eqn:traceson:linearizationGradient} directly from
  \eqref{eqn:traceson:computationHessian}.
  More precisely, utilize
  \begin{equation}
    \begin{aligned}
      H_{f}(\Theta_0)(X,X)
      &=
      \tr{X^{\T}\Theta_0X}
      =
      \vectorize(\Theta_0{X}^{\T})\vectorize(X)
      \\
      &=
      \vectorize(X)^{\T}(I\otimes\Theta_0+(I\otimes{\Theta_0})^{\T})\vectorize(X)
      \\
      &=
      -\frac{1}{2}\vectorize(X)^{\T}\vectorize(\Theta_0^{\T}X+X\Theta_0)
      \text{.}
    \end{aligned}
    \label{}
  \end{equation}
  If $\Theta_0=I$, then the linearization is
  \begin{equation}
    \begin{aligned}
      \dot{X}
      &=
      -X
      \text{,}
    \end{aligned}
  \end{equation}
  which shows that the equilibrium $I$ is locally exponentially stable.
  Now consider the linearization at the other equilibrium points, 
  i.e. $\Theta_0\neq{I}$ and $\Theta_0=\Theta_0^{\T}$.
  Since~$\Theta_0$ is symmetric, $\Theta_0$ is orthonormally diagonalizable,
  i.e. $\Theta_0=\Pi^{\T}D\Pi$ for some diagonal~$D$ and orthonormal~$\Pi$ where
  the columns of~$\Pi$ are eigenvectors of~$\Theta_0$.
  Since~$\Theta_0^{\T}\Theta_0=I$, we get~$D^2=I$ and consequently the
  eigenvalues are~$\pm{1}$.  Since~$\Theta_0\neq{I}$ and~$\det(\Theta_0)=1$, we
  always have an even number of negative eigenvalues with associated
  eigenvectors~$v_1,\ldots,v_k$.
  Set $\overline{U}=v_1v_2^{\T}-v_2v_1^{\T}$
  and $X=\overline{U}\Theta_0\in{T_{\Theta_0}SO(n)}$.  Therefore
  \begin{equation}
    \begin{gathered}
      -\frac{1}{2}(\Theta_0X+X\Theta_0)
      =
      -\frac{1}{2}\left(
      \Theta_0^{\T}\overline{U}\Theta_0+\overline{U}\Theta_0\Theta_0
      \right)
      \\
      =
      -\frac{1}{2}\left(
      \Theta_0(v_1v_2^{\T}-v_2v_1^{\T})\Theta_0
      +
      (v_1v_2^{\T}-v_2v_1^{\T})\Theta_0\Theta_0
      \right)
      \\
      =
      -\frac{1}{2}\left(
      (-v_1v_2^{\T}+v_2v_1^{\T})\Theta_0
      +
      (-v_1v_2^{\T}+v_2v_1^{\T})\Theta_0
      \right)
      =
      \overline{U}\Theta_0
      \text{.}
    \end{gathered}
  \end{equation}
  Therefore, $X=\overline{U}\Theta_0$ is an eigenvector of the operator defined
  by the right hand side of \eqref{eqn:traceson:linearizationGradient}.
  Since the associated eigenvalue is positive (one), the linearization
  \eqref{eqn:traceson:linearizationGradient} is unstable.
  Consequently, the linearization of the gradient flow
  at the equilibria $\Theta_0$ with $\Theta_0\neq{=}I$ and
  $\Theta_0=\Theta_0^{\T}$ is
  unstable, which proves \ref{lem:traceson:7:2}.\\
  \ref{lem:traceson:7:3}
  Denote the flow of $\dot{\Theta}=-\grad{f}(\Theta)$ by 
  $\phi:\mathbb{R}\times{SO(n)}\rightarrow{SO(n)}$ and by ${B}_k$ the
  set 
  \begin{equation}
    {B}_k
    =
    \{x_0\in{SO(n)}\rvert{\lim_{t\rightarrow\infty}\phi(t,x_0)\in\mathcal{F}_k}\}
    \label{}
    \text{,}
  \end{equation}
  i.e. the set of initial conditions that converges to the connected component
  $\mathcal{F}_k$ of the set of critical points $\mathcal{F}$ given in Lemma
  \ref{lem:traceson:1}\ref{lem:traceson:1:3}.
  Because of Proposition \ref{prop:traceson:5}, we are certain that any solution
  of the gradient flow converges to the critical set of $f$, and as a
  consequence,
  $SO(n)=\cup_{k=0}^{\lfloor{\tfrac{n}{2}}\rfloor}{{B}_k}$.
  Then ${B}_0$ is the set of
  initial conditions for which the flow converges to $I$ and
  ${B}=\cup_{k=1}^{\lfloor{\tfrac{n}{2}}\rfloor}{B}_k$ is the
  set of initial conditions for which the flow converges to any of the other
  critical points. In Lemma \ref{lem:traceson:6} we showed that ${B}_k$
  has measure zero and that $SO(n)\setminus{{B}_k}$ is dense in $SO(n)$.
  Since ${B}$ is the union of a finite number of sets of measure zero,
  it has measure zeros, see e.g. \cite[Lemma
  10.1]{Lee2006introductionToSmoothManifolds}. In particular
  ${B}_0=SO(n)\setminus{B}$ is dense.
\end{mypf}

\bibliographystyle{elsarticle-num}
\bibliography{references}
\end{document}